\documentclass[11pt,reqno]{amsart}
\usepackage{amsmath,amsfonts, amssymb, amsthm}
\usepackage{mathrsfs}
\usepackage{graphicx, color,hyperref}

\def\qed{\hfill \rule{4pt}{7pt}}

\def\Z{\mathbb{Z}}
\def\N{\mathbb{N}}

\def\l{\left}
\def\r{\right}
\def\bg{\bigg}
\def\({\bg(}
\def\){\bg)}
\def\t{\text}
\def\f{\frac}

\def\ls{\leqslant}

\def\bi{\binom}

\def\qed{\hfill \rule{4pt}{7pt}}
\newtheorem{thm}{Theorem}[section]

\begin{document}
\hbox{Preprint}

\title[Taylor coefficients and series involving harmonic numbers]
      {Taylor coefficients and series involving harmonic numbers}


\author[Q.-H. Hou]{Qing-Hu Hou}
\address{School of Mathematics, Tianjin
University, Tianjin 300072, People's Republic of China}
\email{{\tt qh\_hou@tju.edu.cn}
\newline\indent
{\it Homepage}: {\tt http://faculty.tju.edu.cn/HouQinghu/en/index.htm}}

\author[Z.-W. Sun]{Zhi-Wei Sun}
\address{Department of Mathematics, Nanjing
University, Nanjing 210093, People's Republic of China}
\email{{\tt zwsun@nju.edu.cn}
\newline\indent
{\it Homepage}: {\tt http://maths.nju.edu.cn/\lower0.5ex\hbox{\~{}}zwsun}}

\keywords{Binomial coefficients, combinatorial identities, harmonic numbers, infinite series.
\newline \indent Supported by the Natural Science Foundation of China (grant no. 11921001 and 12371004, respectively).}

\begin{abstract} During 2022--2023 Z.-W. Sun posed many conjectures on infinite series with summands involving generalized harmonic numbers. Motivated by this, we deduce $58$ series identities involving harmonic numbers, eight of which were previously conjectured by the second author. For example, we obtain that
\[
\sum_{k=1}^{\infty} \frac{(-1)^k}{k^2{2k \choose k}{3k \choose k}}
\left( \frac{7 k-2}{2 k-1} H_{k-1}^{(2)}-\frac{3}{4 k^2} \right) = \frac{\pi^4}{720}.
\]
and
\[
\sum_{k=1}^\infty \frac{1}{k^2 {2k \choose k}^2} \left( \frac{30k-11}{k(2k-1)}  (H_{2k-1}^{(3)} + 2 H_{k-1}^{(3)}) + \frac{27}{8k^4} \right) = 4 \zeta(3)^2,
\]
where $H_n^{(m)}$ denotes $\sum_{0<j\ls n}j^{-m}$.
\end{abstract}
\maketitle

\section{Introduction}

For each $m\in\Z^+=\{1,2,3,\ldots\}$, the $m$th harmonic numbers are those rational numbers
$$H_n^{(m)}=\sum_{0<k\le n}\f1{k^m}\ \ (n\in\N=\{0,1,2,\ldots\}).$$
The numbers $H_n=H_n^{(1)}\ (n\in\N)$ are the usual harmonic numbers.
Z.-W. Sun \cite{Sun22,Sun24, Sun26} formulated many conjectural series identities whose summands involve
generalized harmonic numbers of small orders. Motivated by this,
we confirmed some conjectures of Sun and evaluate more series involving harmonic numbers of small orders.

Pilehroods \cite{PP12} used the WZ method to prove the following identity proposed by the second author in  \cite{Sun11}:
\begin{equation*}\label{27-a}
	\sum_{k=1}^\infty \frac{(-27)^{k-1}(15k-4)}{k^3  {2k \choose k}^2{3k \choose k}} = \sum_{k=1}^\infty \frac{\left( \frac{k}{3} \right)}{k^2},
\end{equation*}
where $\left(-\right)$ denotes the Legendre symbol.
Motivated by their method, we confirm Conjecture 24(i) in \cite{Sun22}. Namely, we have the following theorem.
\begin{thm}\label{bi-27} We have the identity
\begin{equation}\label{eq-bi-27}
\sum_{k=1}^\infty \frac{(-27)^{k}}{k^3 {2k \choose k}^2 {3k \choose k}} \left(  (15k-4) (3 H_{3k-1} - H_{k-1}) - 9 \right) = -\f{4\pi^3}{\sqrt3}.
\end{equation}
\end{thm}

We also have
\begin{thm} \label{bi3}
We have
\begin{equation}\label{eq-16bi3}
	\sum_{k=1}^\infty \frac{16^k}{k^3 {2k \choose k}^3} \left( (6k-2) H_{k-1}^{(2)} + \frac{1}{k} \right) = \frac{\pi^4}{24}.
\end{equation}
\end{thm}

Recently, Wei and Xu \cite{WX2308} confirmed the idenities
\begin{equation*}\sum_{k=0}^\infty\f{\bi{2k}k^5}{(-4096)^k}\l((20k^2+8k+1)H_k^{(3)}+\f 8{2k+1}\r)=\f{64\zeta(3)}{\pi^2}
 \end{equation*}
 and
\begin{equation*}\sum_{k=0}^\infty\f{\bi{2k}k^5}{(-2^{20})^k}\l((820k^2+180k+13)\l(9H_{2k}^{(3)}-H_k^{(3)}\r)+\f {125}{2k+1}\r)=1024\f{\zeta(3)}{\pi^2},
 \end{equation*}
 conjectured by Sun. Our third theorem gives a similar result involving harmonic numbers of order $4$.

\begin{thm}\label{bi5-12}
We have
\begin{equation}\label{eq-bi5-20}
		\sum_{k=0}^\infty \frac{{2k \choose k}^5} {(-2^{20})^k} \left( (820k^2+180k+13) (49 H_{2k}^{(4)} - 3H_k^{(4)}) - \frac{195}{(2k+1)^2} \right)
		= -\frac{896}{45} \pi^2.
\end{equation}
\end{thm}

Chu and Zhang \cite{ChZh} presented several transformation formulas. By taking suitable linear combinations of the Taylor coefficients and using Au's package on multiple zeta functions \cite{Au}, we are able to derive some identities involving harmonic numbers which are similar to some conjectures due to Sun \cite{Sun24, Sun26}.

\cite[Example 84]{ChZh} has the equivalent form:
\begin{equation}\label{Ex84}
\sum_{k=1}^{\infty} \frac{16^k(6 k-1)}{(2 k-1) k^2{2k \choose k}{4k \choose 2k}}  = 8G,
\end{equation}
where $G=\sum_{k=0}^\infty(-1)^k/(2k+1)^2$ is the Catalan constant.
In contrast, we obtain the following result.

\begin{thm}\label{16-b2b4} We have
\begin{equation}\label{eq-16-b2b4}		
\sum_{k=1}^{\infty} \frac{16^k(6 k-1)(4 H_{2 k-1}-H_{k-1})}{(2 k-1) k^2{2k \choose k}{4k \choose 2k}}
		=\pi^3,
\end{equation}
\begin{equation}\label{eqt-1.4b}		
	\sum_{k=1}^{\infty} \frac{16^k \left( (6 k-1)(2 H_{4 k-1}-H_{2k-1}) + \frac{8k}{2k-1} \right) } {(2 k-1) k^2{2k \choose k}{4k \choose 2k}}
	=\pi^3,
\end{equation}
and
\begin{equation}\label{eqt-1.4c}		
	\sum_{k=1}^{\infty} \frac{16^k \left( (6 k-1) H_{k-1}^{(2)} - \frac{8(4k-1)}{(2k-1)^2}  \right) }{ (2 k-1) k^2{2k \choose k}{4k \choose 2k}}
	=-32 \beta(4),
\end{equation}
where 
\[
\beta(4) = \sum_{k=0}^\infty \frac{(-1)^k}{(2k+1)^4}.
\]
\end{thm}

\cite[Example 11]{ChZh} has the equivalent form:
\begin{equation}\label{Ex11}
\sum_{k=1}^{\infty} \frac{30 k-11}{(2 k-1) k^3{2k \choose k}^2}  = 4\zeta(3).
\end{equation}
Wei \cite{Wei2303} deduced some variants of this identity involving harmonic numbers of order $1$ or $2$.
In contrast, we obtain the following result.

\begin{thm}\label{bi2} We have
\begin{equation}\label{eq-bi2}
	\sum_{k=1}^\infty \frac{1}{k^2 {2k \choose k}^2} \left( \frac{30k-11}{k(2k-1)}  (H_{2k-1}^{(3)} + 2 H_{k-1}^{(3)}) + \frac{27}{8k^4} \right) = 4 \zeta(3)^2.
\end{equation}
\end{thm}

\cite[Example 21]{ChZh} has the equivalent form:
\begin{equation}\label{Ex21}
\sum_{k=1}^\infty\f{(-1)^{k-1}(56k^2-32k+5)}{(2k-1)^2k^3\bi{2k}k\bi{3k}k}=4\zeta(3).
\end{equation}
Motivated by this and the spirit of Sun \cite{Sun24, Sun26},  we establish the following theorem.

\begin{thm}\label{Chu21} We have
	\begin{equation}\label{eq-Chu21}
	\sum_{k=1}^\infty \frac{(-1)^{k-1}}{k^3 {2k \choose k} {3k \choose k}}
	\left( \frac{56k^2-32k+5}{(2k-1)^2} H_{2k-1} + \frac{1}{2k} \right)	=  \frac{\pi^4}{20}
\end{equation}
and
\begin{equation}\label{eq-1.12}
\sum_{k=1}^\infty \frac{(-1)^k	\big( (56k^2-32k+5) (H_{3k-1}-\frac{46}{3} H_{2k-1} - H_{k-1}) -12k \big)}{(2k-1)^2 k^3 {2k \choose k} {3k \choose k}} =\frac{7}{10} \pi^4.
\end{equation}
\end{thm}

\cite[Example 24]{ChZh} has the equivalent form:
\begin{equation}\label{Ex24}
\sum_{k=1}^\infty\f{(-1)^k(7k-2)}{(2k-1)k^2\bi{2k}k\bi{3k}k}=-\f{\pi^2}{12}.
\end{equation}
Motivated by this and the spirit of Sun \cite{Sun24, Sun26},  we establish the following theorem.

\begin{thm}\label{Chu24} We have
\begin{equation}\label{eq-Chu24a}
	\sum_{k=1}^\infty \frac{(-1)^k}{k^2 {2k \choose k} {3k \choose k}}
	\left( \frac{7k-2}{2k-1} (2 H_{2k-1} - H_{k-1}) + \frac{1}{k} \right)	= - \frac{3}{2} \zeta(3)
\end{equation}
and
\begin{equation}\label{eq-Chu24b}
	\sum_{k=1}^{\infty} \frac{(-1)^k}{k^2{2k \choose k}{3k \choose k}}\left(\frac{7 k-2}{2 k-1} H_{k-1}^{(2)}-\frac{3}{4 k^2}\right)=\frac{\pi^4}{720}.
\end{equation}
\end{thm}

\cite[Example 1]{ChZh} has the equivalent form:
\begin{equation}\label{Ex1}
\sum_{k=1}^\infty
\frac{(-16)^k (20k^2-8k+1) }{(2k-1)^2 k^3 {2k \choose k} {4k \choose 2k}} = - 14 \zeta(3).
\end{equation}
Motivated by this and the spirit of Sun \cite{Sun24, Sun26},  we establish the following theorem.

\begin{thm}\label{th-1.9}
We have
\begin{equation}\label{eq-th8-a}
\sum_{k=1}^\infty \frac{(-16)^k	\big( (20k^2-8k+1) H_{4k-1}  + \frac{k(28k-5)}{2k-1} \big)}{(2k-1)^2 k^3 {2k \choose k} {4k \choose 2k}} = - \frac{5 \pi^4}{8},
\end{equation}
\begin{equation}\label{eq-th8-b}
	\sum_{k=1}^\infty \frac{(-16)^k	\big( (20k^2-8k+1) H_{2k-1}  + \frac{2k(4k-1)}{2k-1} \big)}{(2k-1)^2 k^3 {2k \choose k} {4k \choose 2k}} = - \frac{\pi^4}{4},
\end{equation}
\begin{equation}\label{eq-th8-c}
	\sum_{k=1}^\infty \frac{(-16)^k	\big( (20k^2-8k+1) (2H_{4k-1}-9H_{2k-1}) -8k \big)}{(2k-1)^2 k^3 {2k \choose k} {4k \choose 2k}} = \pi^4,
\end{equation}
and
\begin{equation}\label{eq-th8-d}
	\sum_{k=1}^\infty \frac{(-16)^k	\big( (20k^2-8k+1) (16H_{2k-1}^{(3)}-3H_{k-1}^{(3)}) - \frac{32k(4k-1)}{(2k-1)^3} \big)}{(2k-1)^2 k^3 {2k \choose k} {4k \choose 2k}} = - 98 \zeta(3)^2.
\end{equation}
\end{thm}

\cite[Example 50]{ChZh} has the equivalent form:
\begin{equation}\label{Ex50}
\sum_{k=1}^\infty\f{256^k(22k-1)}{(2k-1)k^2\bi{3k}k\bi{6k}{3k}}=128G.
\end{equation}
Motivated by this and the spirit of Sun \cite{Sun24, Sun26},  we establish the following theorem.

\begin{thm}\label{th-22}
We have
\begin{equation}\label{eqt-1.9a}
	\sum_{k=1}^\infty \frac{256^k	\big( (22k-1) H_{2k-1} +15 \frac{6k-1}{2k-1} \big)}{(2k-1) k^2 {3k \choose k} {6k \choose 3k}}= 16 \pi^3,
\end{equation}
and
\begin{equation}\label{eqt-1.9b}
	\sum_{k=1}^\infty \frac{256^k	\big( (22k-1) (2H_{6k-1} - H_{3k-1} +H_{k-1}) + \frac{268k-38}{2k-1} \big)}{(2k-1) k^2 {3k \choose k} {6k \choose 3k}}= 48 \pi^3.
\end{equation}
Consequently, we have
\begin{equation}\label{eqt-1.9c}
	\sum_{k=1}^\infty \frac{256^k	\big( (22k-1)\mathcal H(k)-50 \big)}{(2k-1) k^2 {3k \choose k} {6k \choose 3k}}= -16 \pi^3,
\end{equation}
where $\mathcal H(k) = 10H_{6k-1}-5H_{3k-1}-16H_{2k-1}+5H_{k-1}$.

Moreover, we have
\begin{equation}\label{eqt-1.9d}
	\sum_{k=1}^\infty \frac{256^k	\big( (22k-1) H_{2k-1}^{(2)} - 29 \frac{6k-1}{(2k-1)^2} \big)}{(2k-1) k^2 {3k \choose k} {6k \choose 3k}}= -512 \beta(4),
\end{equation}
where 
\[
\beta(4) = \sum_{k=0}^\infty \frac{(-1)^k}{(2k+1)^4}.
\]
\end{thm}

\cite[Example 41]{ChZh} has the equivalent form:
\begin{equation}\label{Ex41}
\sum_{k=1}^\infty\f{(-256)^k(86k^2-21k+2)}{(2k-1)^2k^3\bi{3k}k\bi{6k}{3k}}=-224\zeta(3).
\end{equation}
Motivated by this and the spirit of Sun \cite{Sun24, Sun26},  we establish the following theorem.

\begin{thm}\label{th-1.10}
We have
\begin{equation}\label{eq-4.33}
	\sum_{k=1}^\infty \frac{(-256)^k	\big( (86k^2-21k+2) H_{2k-1} +25k(6k-1)/(2k-1) \big)}{(2k-1)^2 k^3 {3k \choose k} {6k \choose 3k}}= - 8 \pi^4
\end{equation}
and
\begin{equation}\label{eq-4.34}
	\sum_{k=1}^\infty \frac{(-256)^k	\big( (86k^2-21k+2) c_k + 2k(322k-45)/(2k-1) \big)}{(2k-1)^2 k^3 {3k \choose k} {6k \choose 3k}}= - 32 \pi^4,
\end{equation}
where $c_k = 2H_{6k-1}-H_{3k-1}+H_{k-1}$.
Also,
\begin{equation}\label{eq-4.35}
	\sum_{k=1}^\infty \frac{(-256)^k	\big( (86k^2-21k+2) H_{2k-1}^{(2)} - 61k(6k-1)/(2k-1)^2 \big)}{(2k-1)^2 k^3 {3k \choose k} {6k \choose 3k}}= 992 \zeta(5)
\end{equation}
and
\begin{multline}\label{eq-4.36}
	\sum_{k=1}^\infty \frac{(-256)^k	\big( (86k^2-21k+2) (15 H_{2k-1}^{(3)} - 2 H_{k-1}^{(3)}) - 83k(6k-1)/(2k-1)^3 \big)}{(2k-1)^2 k^3 {3k \choose k} {6k \choose 3k}}\\
	= -1568 \zeta(3)^2.
\end{multline}
\end{thm}

\cite[Example 93]{ChZh} has the equivalent form:
\begin{equation}\label{Ex93}
\sum_{k=1}^\infty\f{(-1)^{k-1}(112k^3-8k^2-6k+1)\bi{2k}k^2}{(2k-1)^3k^2\bi{3k}k\bi{6k}{3k}}=\f23\pi^2.
\end{equation}
Motivated by this and the spirit of Sun \cite{Sun24, Sun26},  we establish the following theorem.

\begin{thm}\label{th-1.11}
Let $P(k) = 112k^3-8k^2-6k+1$. Then
\begin{equation}\label{eq-1.11a}
\sum_{k=1}^\infty \frac{(-1)^{k-1} {2k \choose k}^2 \big( P(k) H_{2k-1} + 20k^2(6k-1)/(2k-1) \big)}{(2k-1)^3 k^2 {3k \choose k} {6k \choose 3k}} = 11 \zeta(3)
\end{equation}
and
\begin{equation}\label{eq-1.11b}
	\sum_{k=1}^\infty \frac{(-1)^{k-1} {2k \choose k}^2 \big( P(k) (2 H_{6k-1}-H_{3k-1}+H_{k-1}) + 8k^2 \frac{86k-11}{2k-1} \big)}{(2k-1)^3 k^2 {3k \choose k} {6k \choose 3k}} = 50 \zeta(3).
\end{equation}
\end{thm}

\cite[Example 25]{ChZh} has the equivalent form:
\begin{equation}\label{Ex25}
	\sum_{k=1}^\infty\f{(-16)^{k}(112k^3-116k^2+35k-3)}{k^2(2k-1)(4k-1)(4k-3) \bi{3k}k\bi{6k}{3k}}=- \frac{\pi^2}{4}.
\end{equation}
Motivated by this and the spirit of Sun \cite{Sun24, Sun26}, we establish the following theorem.

\begin{thm}\label{th-1.12}
	Let $P(k) = 112k^3-116k^2+35k-3$. Then
	\begin{equation}\label{eq-1.12a}
		\sum_{k=1}^\infty \frac{(-16)^{k} \big( P(k) H_{2k-1} -\frac{ (4k-1)(6k-1)(8k-3) }{4(2k-1)} \big)}{k^2(2k-1)(4k-1)(4k-3) \bi{3k}k\bi{6k}{3k}} = \frac{\pi^2}{2} \log 2 - \frac{7}{2} \zeta(3)
	\end{equation}
	and
	\begin{multline}\label{eq-1.12b}
		\sum_{k=1}^\infty \frac{(-16)^{k} \big( P(k) (2H_{6k-1}-H_{3k-1}-3H_{k-1}) -64k^2+46k -17/2 \big)}{k^2(2k-1)(4k-1)(4k-3) \bi{3k}k\bi{6k}{3k}} \\
		= 7 \zeta(3) -2\pi^2 \log 2.
	\end{multline}
\end{thm}

\cite[Example 72]{ChZh} has the equivalent form:
\begin{equation}\label{Ex72}
	\sum_{k=1}^\infty\f{(3k-1)4^{k}}{(2k-1)k^3 {2k \choose k}^2}= \frac{7}{4} \zeta(3).
\end{equation}
Motivated by this and the spirit of Sun \cite{Sun24, Sun26},  we establish the following theorem.
\begin{thm}\label{th-1.13}
We have
\begin{equation}\label{eq-1.13a}
\sum_{k=1}^\infty \frac{4^k}{(2k-1)k^3{2k \choose k}^2} \left( (3k-1)(2H_{2k-1}-3H_{k-1}) + \frac{2k}{2k-1}\right)
= \frac{\pi^4}{16},
\end{equation}
\begin{equation}\label{eq-1.13b}
	\sum_{k=1}^\infty \frac{(3k-1)4^k}{(2k-1)k^3{2k \choose k}^2} \left( 4H^{(2)}_{2k-1} - 5H_{k-1}^{(2)} \right)
	= \frac{31}{4} \zeta(5),
\end{equation}
and
\begin{equation}\label{eq-1.13c}
	\sum_{k=1}^\infty \frac{(3k-1)4^k}{(2k-1)k^3{2k \choose k}^2} \left(8H^{(3)}_{2k-1} +7H_{k-1}^{(3)} \right)
= \frac{49}{4} \zeta(3)^2.
\end{equation}
\end{thm}

\cite[Example 47]{ChZh} has the equivalent form:
\begin{equation}\label{Ex47}
	\sum_{k=1}^\infty\f{(22k^2-17k+3) 16^k {4k \choose 2k}}{k (4k-1)(4k-3) {3k \choose k}{6k \choose 3k}}= 2 \pi.
\end{equation}
Motivated by this and the spirit of Sun \cite{Sun24, Sun26},  we establish the following theorem.
\begin{thm}\label{th-1.14}
	Let $Q(k) = (2k-1)(11k-3) = 22k^2-17k+3$. 
	We have
	\begin{equation}\label{eq-1.14a}
		\sum_{k=1}^\infty \frac{16^k {4k \choose 2k} \big( Q(k) H_{k-1} + (6k-1)(4k-3)/(2k-1) \big) }{k (4k-1)(4k-3){3k \choose k} {6k \choose 3k}} 
		= 16G - 4 \pi \log 2,
	\end{equation}
	\begin{equation}\label{eq-1.14b}
		\sum_{k=1}^\infty \frac{16^k {4k \choose 2k} \big( Q(k) H_{2k-1} + (6k-1)(k-3)/(6k-3) \big) }{k (4k-1)(4k-3){3k \choose k} {6k \choose 3k}} 
		= \frac{4}{3} (8G - \pi \log 2),
	\end{equation}
	and
	\begin{multline}\label{eq-1.14c}
	\sum_{k=1}^\infty \frac{16^k {4k \choose 2k} \big( Q(k) (2H_{6k-1}-H_{3k-1})+ 4(18k^2-20k+3)/(6k-3) \big) }{k (4k-1)(4k-3){3k \choose k} {6k \choose 3k}} 
	\\ =  \frac{4}{3} (16G + \pi \log 2).
	\end{multline}
\end{thm}

\cite[Example 14]{ChZh} has the equivalent form:
\begin{equation}\label{Ex14}
	\sum_{k=1}^\infty\f{(60k^2-43k+8)  {4k \choose 2k}} {k^3(4k-1) {2k \choose k}^4}= \f{\pi^2}{3}.
\end{equation}
Motivated by this and the spirit of Sun \cite{Sun24, Sun26},  we establish the following theorem.
\begin{thm}\label{th-1.15}
	Let $Q(k) = 60k^2-43k+8$. 
	We have
	\begin{equation}\label{eq-1.15a}
		\sum_{k=1}^\infty \frac{{4k \choose 2k} \big( Q(k) (3H_{2k-1}-2H_{k-1}) -34k+25/2 \big) }{(4k-1)k^3{2k \choose k}^4} 
		= 6 \zeta(3),
	\end{equation}
	\begin{equation}\label{eq-1.15b}
		\sum_{k=1}^\infty \frac{{4k \choose 2k} \big( Q(k) H_{4k-1} -(304k^2-216k+39)/(4(4k-1)) \big) }{(4k-1)k^3{2k \choose k}^4} 
= 4 \zeta(3),
	\end{equation}
\begin{equation}
\label{eq-1.15c}
\sum_{k=1}^\infty \frac{{4k \choose 2k} \big( Q(k) H_{k-1}^{(2)} +(44k-15)/(4k) \big) }{(4k-1)k^3{2k \choose k}^4} 
= \frac{\pi^4}{90},
\end{equation}
and
\begin{equation}\label{eq-1.15d}
	\sum_{k=1}^\infty \frac{{4k \choose 2k} \big( Q(k) H_{2k-1}^{(2)} +(52k-21)/(4k) \big) }{(4k-1)k^3{2k \choose k}^4} 
	= \frac{2\pi^4}{45}.
\end{equation}
\end{thm}

\cite[Example 34]{ChZh} has the equivalent form:
\[
\sum_{k=1}^\infty  \frac{4096^k (828k^{3}-756k^{2}+199k-15)}{k^3(3k-1)(3k-2){3k \choose k} {6k \choose 3k}^2}  = 48 \pi^2.
\]
Motivated by this and the spirit of Sun \cite{Sun24, Sun26},  we establish the following theorem.
\begin{thm}\label{th-16}
Let $P(k) = 828k^{3}-756k^{2}+199k-15$. We have
\begin{multline} \label{eq-16a}
\sum_{k=1}^{\infty} \frac{4096^{k} (5P(k)(H_{2k-1}-H_{k-1})-1656k^{2}+576k-53} {k^{3}(3k-1)(3k-2){3k \choose k} {6k \choose 3k}^2}  \\
= 96 (3\pi^{2}\log2-14\zeta(3))
\end{multline}
and
\begin{multline} \label{eq-16b}
\sum_{k=1}^{\infty} \frac{4096^{k} (5P(k)(4H_{6k-1}-3H_{3k-1}-H_{k-1})-Q(k)/((3k-1)(3k-2))} {k^{3}(3k-1)(3k-2){3k \choose k} {6k \choose 3k}^2}  \\
= 192 (7 \zeta(3) + 6 \pi^{2}\log2),
\end{multline}
where $Q(k) = 22356k^{4}-38232k^{3}+22617k^{2}-5676k+509$.
\end{thm}

\cite[Example 52]{ChZh} has the equivalent form:
\[
\sum_{k=1}^\infty {\frac { \left( 74{k}^{2}-47k+8 \right) {3k\choose k}{6k
			\choose 3k}}{ \left( 6k-1 \right) {k}^{3}  {2k\choose k}
	 ^{5}}}
 = 2 \pi^2.
\]
Motivated by this and the spirit of Sun \cite{Sun24, Sun26},  we establish the following theorem.
\begin{thm}\label{th-17} 
We have
\begin{equation}\label{eq-17a}
\sum_{k=1}^\infty {\frac { {3k\choose k}{6k
			\choose 3k}}{ \left( 6k-1 \right) {k}^{3}  {2k\choose k}^{5}}}
		\big( 2  \left( 74{k}^{2}-47k+8 \right) H_{2k-1} - 34k + 11 \big)
= 36 \zeta(3).		 
\end{equation}
\begin{multline} \label{eq-17b}
\sum_{k=1}^\infty {\frac { {3k\choose k}{6k
			\choose 3k}}{ \left( 6k-1 \right) {k}^{3}  {2k\choose k}^{5}}} 
\left( (74 k^2-47 k+8) (200 H_{2k-1}^{(2)} - 137 H_{k-1}^{(2)}) - \frac{168}{k}  \right) \\
= \frac{526}{15} \pi^4.	
\end{multline}
\begin{multline} \label{eq-17c}
	\sum_{k=1}^\infty {\frac { {3k\choose k}{6k
				\choose 3k}}{ \left( 6k-1 \right) {k}^{3}  {2k\choose k}^{5}}} 
	\Bigg( (74 k^2-47 k+8) ( 7H_{2k-1}^{(3)} - 3 H_{k-1}^{(3)}) \\
	+ \frac{1228k^2-648k+81}{8k^2(2k-1)}  \Bigg) 
	= 8\pi^2 \zeta(3) +60 \zeta(5).	
\end{multline}

\end{thm}

\cite[Example 4]{ChZh} has the equivalent form:
\[
\sum_{k=0}^\infty \frac { \left(20k^2+44k+25 \right) {2k\choose k}}{ (2k+1)^4 (2k+3)^4 (-16)^k}
= \frac{\pi^2}{32}.
\]
Motivated by this and the spirit of Sun \cite{Sun24, Sun26},  we establish the following theorem.
\begin{thm}\label{th-18} 
We have
\begin{multline}
\label{eqt-18}
\sum_{k=0}^\infty \frac { {2k\choose k} \left( \left( 20k^2+44k+25 \right) (H_{2k}-H_k) -2\,\frac{120\,{k}^{3}+388\,{k}^{2}+434\,k+167}{ \left( 2\,k+1
				\right)  \left( 2\,k+3 \right) }
		 \right)}{ (2k+1)^4 (2k+3)^4 (-16)^k} \\
= \frac{1}{16}\,{\pi }^{2} \ln  \left( 2 \right) - \frac{3}{4} - {\frac {7\,\zeta  \left( 3
		\right) }{8}}
\end{multline}
\end{thm}

\cite[Example 101]{ChZh} has the equivalent form:
\[
\sum_{k=1}^\infty {\frac {{64}^{k} (368 k^3-396 k^2+116 k-9)}{{k}^{2} \left( 2\,k-1 \right)  \left( 4\,k-1
		\right)  \left( 4\,k-3 \right) {3\,k\choose k}{6\,k\choose 3\,k}}}
= 32 G.
\]
Motivated by this and the spirit of Sun \cite{Sun24, Sun26},  we establish the following theorem.
\begin{thm}\label{th-19} 
We have
\begin{equation}\label{eqt-19}
\sum_{k=1}^\infty \frac { {64}^{k}   \Big(   (368 k^3-396 k^2+116 k-9) {\mathcal H}(k)  + \frac {240 k^3-252 k^2+12 k+11}{2\,k-1} \Big) } { {k}^{2} \left( 2\,k-1 \right)  \left( 4\,k-1\right)  \left( 4\,k-3 \right) {3\,k\choose k} {6\,k\choose 3\,k} }
= 4 \pi^3,
\end{equation}
where ${\mathcal H}(k) = 2 H_{6k-1} - H_{3k-1} + H_{2k-1} - H_{k-1}$.
\end{thm}


\cite[Example 85]{ChZh} has the equivalent form:
\[
\sum_{k=1}^{\infty} \frac{(40k^2 - 24k + 3)(-256)^k}{k^3(2k-1)\binom{2k}{k}\binom{4k}{2k}^2}
= -64 G.
\]
Motivated by this and the spirit of Sun \cite{Sun24, Sun26},  we establish the following theorem.
\begin{thm}\label{th-20} 
	We have
	\begin{equation}\label{eqt-20a}
\sum_{k=1}^{\infty} \frac{(-256)^k \left( (40k^2 - 24k + 3)(2H_{4k-1} - H_{2k-1}) - 12k + 4) \right)}{k^3(2k-1)\binom{2k}{k}\binom{4k}{2k}^2} = -4 \pi^3,
	\end{equation}
\begin{equation}\label{eqt-20b}
\sum_{k=1}^{\infty} \frac{(-256)^k \left( (40k^2 - 24k + 3)(5H_{2k-1} - 3H_{k-1}) - \frac{88k^2 - 48k + 7}{2k-1} \right)}{k^3(2k-1)\binom{2k}{k}\binom{4k}{2k}^2} = -4 \pi^3,
\end{equation}
\begin{equation}
\sum_{k=1}^{\infty} \frac{(-256)^k \left( (40k^2 - 24k + 3)(4H_{4k-1}^{(2)} - H_{2k-1}^{(2)} - H_{k-1}) - 4 \right)}{k^3(2k-1)\binom{2k}{k}\binom{4k}{2k}^2}
= -256 \beta(4).
\end{equation}

\end{thm}

\cite[Example 63]{ChZh} has the equivalent form:
\[
\sum_{k=1}^{\infty} \frac{(-1)^{k-1}(410k^2 - 197k + 24)}{k^3(2k-1)\binom{2k}{k}\binom{4k}{2k}^2} = \frac{\pi^2}{3}.
\]
Motivated by this and the spirit of Sun \cite{Sun24, Sun26},  we establish the following theorem.

\begin{thm}\label{th-21} 
Let $P(k) = 410k^2-197k +24$. 
We have
\begin{equation}\label{eq-21a}
\sum_{k=1}^{\infty} \frac{(-1)^k \left( P(k) (H_{2k-1} - H_{k-1}) - \frac{1148k^2 - 572k + 7}{4k - 2} \right)}{k^3(2k-1)\binom{2k}{k}\binom{4k}{2k}^2} = \zeta(3),	
\end{equation}
\begin{equation}\label{eq-21b}
\sum_{k=1}^{\infty} \frac{(-1)^k \left( (P(k)(H_{4k-1} - H_{2k-1}) + \frac{164k^2 + 124k - 39}{8k - 4} \right)}{k^3(2k-1)\binom{2k}{k}\binom{k4}{2k}^2} = -3\zeta(3),
\end{equation}
\begin{equation}\label{eq-21c}
\sum_{k=1}^{\infty} \frac{(-1)^{k-1}( (P(k)(4H_{4k-1}^{(2)} - H_{2k-1}^{(2)} - 3H_{k-1}^{(2)}) - 124+ \frac{6}{k})}{k^3(2k-1)\binom{2k}{k}\binom{4k}{2k}^2} = \frac{2}{15}\pi^4,
\end{equation}
\begin{equation}\label{eq-21d}
\sum_{k=1}^{\infty} \frac{(-1)^{k-1}((P(k)(12H_{4k-1}^{(2)} - 27H_{2k-1}^{(2)} - 17H_{k-1}^{(2)}) + 48F(k))}{k^3(2k-1)\binom{2k}{k}\binom{4k}{2k}^2} = \frac{26}{45}\pi^4,
\end{equation}
where $F(k)=(228k^2-100k+9)/(2k -1)^2$.
	
\end{thm}

\cite[Example 12]{ChZh} has the equivalent form:
\[
\sum_{k=1}^{\infty} \frac{256^k (60k^2-26k+3)}{(2k-1) k^4 {2k \choose k}^2 {4k \choose 2k}^2} = 56 \zeta(3).
\]
Motivated by this and the spirit of Sun \cite{Sun24, Sun26},  we establish the following theorem, which confirms Conjecture 5.29 of \cite{Sun26}.

\begin{thm}\label{th-22n} 
	Let $P(k) = 60k^2-26k+3$. 
	We have
	\begin{equation}\label{eq-22a}
		\sum_{k=1}^{\infty} \frac{256^k \left(P(k)(5H_{2k-1} - 2H_{k-1}) - \frac{9(4k-1)^2}{2k-1}\right)}{(2k-1)k^4 \binom{2k}{k}^2 \binom{4k}{2k}^2} = 2\pi^4 ,	
	\end{equation}
	\begin{equation}\label{eq-22b}
		\sum_{k=1}^{\infty} \frac{256^k \left(P(k)(2H_{4k-1} - H_{2k-1}) - \frac{4(2k^2-5k+1)}{2k-1}\right)}{(2k-1)k^4 \binom{2k}{k}^2 \binom{4k}{2k}^2} = 2\pi^4,
	\end{equation}
	\begin{equation}\label{eq-22c}
		\sum_{k=1}^{\infty} \frac{256^k \left(P(k)H_{k-1}^{(2)} - \frac{136k^3 - 76k^2 + 14k - 1}{k(2k-1)^2}\right)}{(2k-1)k^4 \binom{2k}{k}^2 \binom{4k}{2k}^2} = -124\zeta(5) ,
	\end{equation}
	\begin{equation}\label{eq-22d}
		\sum_{k=1}^{\infty} \frac{256^k \left(P(k)\left(H_{4k-1}^{(2)} - \frac{H_{2k-1}^{(2)}}{4}\right) - \frac{Q(k)}{k(2k-1)^2}\right)}{(2k-1)k^4 \binom{2k}{k}^2 \binom{4k}{2k}^2} = 0 ,
	\end{equation}
	where $Q(k) = 92k^3 - 64k^2 + 15k - 1$.
\end{thm}

\cite[Example 118]{ChZh} has the equivalent form:
\[
\sum_{k=1}^{\infty} \frac{364k^2 - 227k + 36}{(2k-1) k^4 {2k \choose k}^2 {3k \choose k}^2} = 4 \zeta(3).
\]
Motivated by this and the spirit of Sun \cite{Sun24, Sun26},  we establish the following theorem, which confirm Conjecture 5.31 of \cite{Sun26}.

\begin{thm}\label{th-23} 
Let  $P(k) = 364k^2 - 227k + 36$. 
We have
	\begin{equation}\label{eq-23a}
		\sum_{k=1}^{\infty} \frac{P(k)(3H_{2k-1} - 2H_{k-1}) - \frac{1276k^2 - 844k + 139}{4k-2}}{(2k-1)k^4 \binom{2k}{k}^2 \binom{3k}{k}^2} = \frac{\pi^4}{15}, 
	\end{equation}
	\begin{equation}\label{eq-23b}
		\sum_{k=1}^{\infty} \frac{P(k)(H_{3k-1} - H_{k-1}) - \frac{728k^2 - 728k + 155}{12k-6}}{(2k-1)k^4 \binom{2k}{k}^2 \binom{3k}{k}^2} = \frac{\pi^4}{15},
	\end{equation}
	\begin{equation}\label{eq-23c}
		\sum_{k=1}^{\infty} \frac{P(k)\left(99H_{3k-1}^{(2)} - 757H_{k-1}^{(2)}\right) + \frac{18Q(k)}{(2k-1)^2}}{(2k-1)k^4 \binom{2k}{k}^2 \binom{3k}{k}^2} = 1316\zeta(5),
	\end{equation}
	where $Q(k) = 2952k^2 - 1572k + 163.$.
\end{thm}

Our method to prove the identities in our theorems is to extract the multivariate Taylor coefficients of a hypergeometric identity. We will utilize two types of identities.

In Section 2, we focus on identities constructed through the WZ method.  Recall that a {\it WZ pair} $(F(n,k), G(n,k))$ refers to a pair of hypergeometric functions that satisfy
\[
F(n+1,k) - F(n,k) = F(n,k+1) - F(n,k), \quad \forall\, n,k \in \mathbb{N}.
\]
It has been demonstrated in \cite{MZ96} that
\begin{equation}\label{WZ-o}
\sum_{k=0}^\infty F(0,k) - \lim_{n \to \infty} \sum_{k=0}^n F(n,k) = \sum_{n=0}^\infty G(n,0) - \lim_{k \to \infty} \sum_{n=0}^k G(n,k).
\end{equation}
In most cases, both limits vanish and $F(n,k)$ and $G(n,k)$ involve additional parameters. By comparing the Taylor expansions of these parameters, we obtain identities that involve harmonic numbers. Following this approach, Au \cite{Au, Au1} confirmed several conjectures made by the second author. We will prove Theorems~\ref{bi-27} and \ref{bi3} in this way. It is worth noting that if $(F(n,k), G(n,k))$ is a WZ pair, then so is $(F(n+a,k+b), G(n+a,k+b))$. Consequently, any WZ pair can be generalized by introducing two extra parameters $a$ and $b$.

In Section 3, we focus on identities established from summation formulas or transformation formulas of hypergeometric series. By applying differential operations, one can derive identities that involve harmonic numbers. This method has been explored by Chu and his coauthor \cite{Chu97, CL23} and Wei \cite{Wei2303, Wei2306}.

We emphasize that symbolic computation plays an important role in our work. We use the {\tt Maple} command {\tt coeftayl} to extract the Taylor coefficients,  the {\tt Maple} package {\tt APCI} (which is available at the first author's homepage) to compute the WZ pair $(F(n,k), G(n,k))$ once $F(n,k)$ is given, and the {\tt Mathematica} package {\tt MultipleZetaValue} developed by Au \cite{Au} to evaluate series involving rational functions and harmonic numbers. We also search for proper linear combinations of the coefficients by solving linear equations.

Let us first introduce some notations.
We use the standard notation of raising factorial
\[
(a)_n = \begin{cases}
	a (a+1) \cdots (a+n-1), &\t{if}\ n>0, \\
	1, &\t{if}\ n=0.
\end{cases}
\]
We will use $[x^i y^j]f(x,y)$ to denote the Taylor coefficient of $x^i y^j$ in $f(x,y)$. If there is no confusion, we write notions as $a[x^2]+b[xy]$ to denote the linear combinations of these coefficients.

\section{Identities from WZ pairs}
\noindent
{\it Proof of Theorem~\ref{bi-27}. }
Let $(F(n,k),G(n,k))$ be the WZ pair given by Pilehroods in \cite{PP12}:
\[
F(n,k) = (-1)^n (n+2k+1) H(n,k),
\]
\[
G(n,k) = (-1)^n \frac{18k^2 + 54kn + 45n^2 + 36k + 63n + 22}{18(2n+1)}  H(n,k),
\]
where
\[
H(n,k) = \frac{(1/3)_k^2 (2/3)_k^2 (2/3)_n (4/3)_n n!^2}{(1/3)_{n+k+1}^2 (2/3)_{n+k+1}^2 {2n \choose n}}.
\]
For the WZ-pair $(F(n,k+b), G(n,k+b))$, we have
\[
\sum_{k=0}^\infty F(0,k+b) = \sum_{n=0}^\infty G(n,b).
\]
The Taylor coefficient $[b^1]$ at $b=0$ of the right hand side is
\[
2\sum_{n=1}^\infty \frac{(-27)^{n}}{n^3 {3n \choose n} {2n \choose n}^2} \left(  (15n-4) (3 H_{3n-1} - H_{n-1}) - 9 \right).
\]
While the Taylor coefficient $[b^1]$ at $b=0$ of the left hand side is
\[
162 \sum_{k=1}^\infty \left( \frac{1}{(3k-1)^3} - \frac{1}{(3k-2)^3} \right)=-162\sum_{k=1}^\infty\f{(\f k3)}{k^3}.
\]
It is known that
$$\sum_{k=1}^\infty\f{(\f k3)}{k^3}=\f{4\pi^3}{81\sqrt3}.$$
Therefore we obtain the desired \eqref{bi-27}. \qed

\noindent
{\it Proof of Theorem~\ref{bi3}. } Let
\begin{align}
	& H(a,b,c,d,k) = (-1)^{a+b}  \frac{\Gamma(a+k+1) \Gamma(b+k+1)}{\Gamma(-a-b+c+d+1))} \nonumber \\
	& \quad \times \frac{\Gamma(-a+c+1) \Gamma(-a+d+1) \Gamma(-b+c+1) \Gamma(-b+d+1)}{\Gamma(c+k+2) \Gamma(d+k+2) } . \label{H}
\end{align}

Au \cite{Au} showed that for
\[
F(n,k) = H(a,b,c+n,d+n,k),
\]
there exists a hypergeometric term $G(n,k)$ such that $(F(n,k), G(n,k))$ forms a WZ pair. One can compute the explicit $G(n,k)$ by the {\tt Maple} package {\tt APCI}. By the WZ pair
\[
\tilde{F}(n,k)=F(n,k)/F(0,0), \quad \tilde{G}(n,k)=G(n,k)/F(0,0),
\]
we derive from \eqref{WZ-o} that
\begin{multline}\label{eq-1.3}
	\sum_{k=0}^\infty \frac{(a+1)_k (b+1)_k}{(c+2)_k (d+2)_k} \\
	= \sum_{n=0}^\infty p_2(n)  \frac{ (-a+c+1)_n (-a+d+1)_n (-b+c+1)_n  (-b+d+1)_n}{(c+2)_n (d+2)_n (-a-b+c+d+1)_{2n+2}},
\end{multline}
where $p_2(n)$ is a polynomial in $n$ with parameters $a,b,c,d$.

Now consider the Taylor expansion in $a,b,c,d$ at $a=b=c=d=-1/2$.
For the left hand side of \eqref{eq-1.3}, we have
\[
[ab]-2[b^2] = \sum_{k=0}^\infty \frac{4H_{2k}^{(2)} - H_k^{(2)}} {(2k+1)^2} .
\]
By {\tt MultipleZetaValues}, we see that this sum is equal to $\pi^4/96$.

For the right hand side of \eqref{eq-1.3}, we have
\[
[ab]-2[b^2] = \frac{1}{4} \sum_{n=1}^\infty  \frac{2(3n-1) H_{n-1}^{(2)} +1/n}{n^3 {2n \choose n}^3} 16^n.
\]
Therefore, \eqref{eq-16bi3} follows immediately.

In view of the above, we have completed the proof of Theorem \ref{bi3}. \qed

\section{Hypergeometric identities}
In this section, we use some hypergeometric identities to show the rest theorems.

\noindent {\it Proof of Theorem~\ref{bi5-12}.}
We quote  Identity (3.6) from \cite{WX2308}:
\begin{multline*}
	\sum_{k=0}^{\infty}  (-1)^k \frac{(a)_k(b)_k(c-b)_k(d-c)_k(1-b)_k(1+b-c)_k(1+c-d)_k}{(1)_{2 k}(1+a-b)_{2 k}(1+a+b-c)_{2 k}(1+a+c-d)_{2 k}} \\
	 \times \frac{(1+a-c)_k(1+a-b+c-d)_k(1+a+b-d)_k}{(2+a-d)_{2 k}} \theta_k(a, b, c, d) \\
	=  \frac{\Gamma(1+a-b) \Gamma(1+a+b-c) \Gamma(1+a+c-d) \Gamma(2+a-d)}{\Gamma(1+a) \Gamma(1+a-c) \Gamma(1+a-b+c-d) \Gamma(1+a+b-d)},
\end{multline*}
where
\begin{multline*}
	 \theta_k(a, b, c, d) 	 = \frac{2 k(1+2 a-d+3 k)}{a} 	 +\frac{(a+k)(1+a-c+k)}{a(1+a-b+2 k)} \\
	 \quad \times \frac{(1+a-b+c-d+k)(1+a+b-d+k)(2+a-b+c-d+3 k)}{(1+a+c-d+2 k)(2+a-d+2 k)} \\
	 +\frac{(a+k)(c-b+k)(1-b+k)(1+c-d+k)}{a(1+2 k)(1+a-b+2 k)(1+a+b-c+2 k)} \\
	 \quad \times \frac{(1+a-c+k)(1+a-b+c-d+k)(1+a+b-d+k)}{(1+a+c-d+2 k)(2+a-d+2 k)} .
\end{multline*}
Let us consider the Taylor expansion at $(a,b,c,d) = (1/2,1/2,1,3/2)$.
By investigating all Taylor coefficients of order $4$, we find the following linear combination
\[
16 [d^4]+12[cd^3]+4[c^2 d^2]+8[ad^3]+4[acd^2]+4[a^2d^2]+2[a^2cd]+2[a^3d]+[a^4].
\]
On the left hand side, the coefficient is
\[
\frac{1}{64} \sum_{k=0}^\infty \frac{{2k \choose k}^5} {(-2^{20})^k} \left( (820k^2+180k+13) (49 H_{2k}^{(4) } - 3H_k^{(4)} + 16) - \frac{195}{(2k+1)^2} \right).
\]
While on the right hand side, the coefficient is
\[
-\frac{14}{45} \pi^2 +  \frac{32}{\pi^2} .
\]
Notice that Guillera have shown in \cite{Gui02} that
\[
\sum_{k=0}^\infty \frac{{2k \choose k}^5} {(-2^{20})^k} (820k^2+180k+13)  = \frac{128}{\pi^2}.
\]
Therefore, \eqref{eq-bi5-20} follows immediately.

In view of the above, the proof of Theorem \ref{bi5-12} is now complete. \qed

\noindent {\it Proof of Theorem~\ref{16-b2b4}.}
Setting $e \to -\infty$ in \cite[Theorem 9]{ChZh}, we get
\begin{align}
	& \sum_{k=0}^\infty (-1)^k (a+2k) \frac{(b)_k (c)_k (d)_k}{(1+a-b)_k (1+a-c)_k (1+a-d)_k} \nonumber \\
	= &\sum_{k=0}^{\infty} \frac{ \alpha_k(a; b, c, d)}{(1+a-b)_{2k}} \cdot
	\frac{(c)_k (d)_k (1+a-b-c)_k (1+a-b-d)_k}{(1+a-c)_k(1+a-d)_k} ,\label{eq-9}
\end{align}
where
{\small \[
		 \alpha_k(a; b, c, d)=\frac {{a}^{2}-ab+4\,ak-2\,bk-cd-ck-kd+3\,{k}^{2}+a+2\,k}{1+a-b+2\,k}
	\]}
Let us consider the linear combination $2[a]+2[b]+[d]$ of Taylor coefficients at $(a,b,c,d)=(1,1/2,1/2,1/2)$.
The coefficient of the right hand side of \eqref{eq-9} is
\[
-\frac{1}{8} \sum_{k=1}^{\infty} \frac{16^k(6 k-1)}{(2 k-1) k^2{2k \choose k}{4k \choose 2k}} (4H_{2k-1} - H_{k-1} - 6).
\]
While the coefficient of the right hand side is
{\small
\[
2\, \sum_{k=0}^\infty   \left( -1 \right) ^{k} {\frac { \left( 6\,k+1 \right)}{ \left( 1+2
		\,k \right) ^{3}}}
\]
}\\
By {\tt MultipleZetaValues}, this sum reduces to $6G-\frac{\pi^3}{8}$. Recall the known identity \eqref{Ex84}. Therefore, \eqref{eq-16-b2b4} follows immediately. 

In a similar way, we can prove \eqref{eqt-1.4b} and \eqref{eqt-1.4c} by the combinations
\[
-16[d]-8[b]-16[a]
\]
and
\[
-16[d^2]+8[cd]. \tag*{\qed}
\]

\noindent {\it Proof of Theorem~\ref{bi2}.}
Recall that \cite[Theorem 14]{ChZh}
\begin{align}
& \sum_{k=0}^\infty (a+2k) \frac{(b)_k (c)_k (d)_k (e)_k}{(1+a-b)_k (1+a-c)_k (1+a-d)_k (1+a-e)_k} \nonumber \\
= &\sum_{k=0}^{\infty} \beta_k(a ; b, c, d, e) \frac{(1+a-b-d)_{2 k}}{(1+a-c_k (1+a-e)_k} \nonumber \\
 &  \times \frac{(c)_k (e)_k (1+a-b-c)_k (1+a-b-e)_k (1+a-c-d)_k (1+a-d-e)_k}{(1+a-b)_{2k} (1+a-d)_{2k} (1+2 a-b-c-d-e)_{2 k}} ,\label{eq-14}
\end{align}
where
{\small \[
\begin{aligned}
& \beta_k(a ; b, c, d, e)=\frac{(1+2 a-b-c-d+3 k)(a-e+k)}{1+2 a-b-c-d-e+2 k} + (e+k)(1+a-b-c+k) \\
	& \qquad  \times \frac{(1+a-b-d+2 k)(1+a-c-d+k)(2+2 a-b-d-e+3 k)}{(1+a-b+2 k)(1+a-d+2 k)(1+2 a-b-c-d-e+2 k)_2}.
\end{aligned}
\]}
Let us consider the Taylor expansion at $(a,b,c,d,e)=(2,1,1,1,1)$.
By investigating all Taylor coefficients of order $3$, we find the following linear combination
\begin{multline*}
[a^3]+[a^2d]+2[a^2e]+[abd]+[ace]-[ad^2]+2[ade]+2[ae^2]+[bd^2]+2[bde]\\
+[cde]+2[ce^2]-2[d^3]-2[d^2e]+2[de^2]+2[e^3].
\end{multline*}
The coefficient of the right hand side of \eqref{eq-14} is
\[
\frac{1}{16}\sum_{k=1}^\infty \frac{1}{k^2 {2k \choose k}^2} \left( \frac{8(30k-11)}{k(2k-1)}  (H_{2k-1}^{(3)} + 2 H_{k-1}^{(3)}) + \frac{27}{k^4} \right).
\]
While the coefficient of the right hand side is
{\small
\begin{multline*}
\sum_{k=0}^\infty \left(
\frac{8 H_k^3}{3 (1+k)^3} 
+ \frac{2 (2 k+1) H_k^2} {(1+k)^4} 
- \frac{4  H_k^{(2)} H_k}{(1+k)^3}
+ \frac{2 k (2 k+3) H_k}{(1+k)^5} \right. \\
\left. - \frac{(2k+1)H_k^{(2)}}{(1+k)^4} 
+ \frac{4 H_k^{(3)}}{3 (1+k)^3} 
- \frac{k (k+3)}{(1+k)^6}
\right).
\end{multline*}}\\
By {\tt MultipleZetaValues}, this sum reduces to $2 \zeta(3)^2$. Therefore, \eqref{eq-bi2} follows immediately. \qed

\noindent {\it Proof of Theorem~\ref{Chu21}.}
Recall that \cite[Theorem 18]{ChZh}
\begin{align}
	& \sum_{k=0}^\infty (a+2k) \frac{(b)_k (c)_k (d)_k (e)_k}{(1+a-b)_k (1+a-c)_k (1+a-d)_k (1+a-e)_k} \nonumber \\
	= &\sum_{k=0}^{\infty} \gamma_k(a ; b, c, d, e) \frac{(-1)^k}{(1+a-b)_{3 k}}  \frac{(c)_k (d)_k (1+a-b-e)_k (1+a-c-d)_k}{(1+a-e)_k} \nonumber \\
	& \times \frac{ (e)_{2k} (1+a-b-c)_{2k} (1+a-b-d)_{2k}}{ (1+a-c, 1+a-d)_{2k} (1+2 a-b-c-d-e)_{2 k}} ,\label{eq-18}
\end{align}
where
{\small \[
	\begin{aligned}
		& \gamma_k(a ; b, c, d, e)= \frac{(1+2 a-b-c-d+4 k)(a-e+k)}{1+2 a-b-c-d-e+2 k} + \frac{ (1+a-b-c+2 k) }{(1+a-b+3 k) } \nonumber \\
			& \times \frac{(e+2 k)(1+a-b-d+2 k) (1+a-c-d+k) (2+2 a-b-d-e+3 k)}{(1+a-d+2 k) (1+2 a-b-c-d-e+2 k)_2} \nonumber \\
			& + \frac{(c+k) (e+2 k) (1+a-b-c+2 k)(1+a-b-e+k) (1+a-c-d+k) }{(1+a-c+2 k) (1+a-d+2 k) } \nonumber \\
			& \times \frac{ (1+a-b-d+2 k)_2}{(1+a-b+3 k)_2 (1+2 a-b-c-d-e+2 k)_2}
		\end{aligned}.
	\]}

To prove the first identity, we consider Taylor expansion at $(a,b,c,d,e)=(2,1,1,1,1)$ and the linear combination
\[
2 [a] + 2 [b] + [e].
\]
The value for the right hand side of \eqref{eq-18} is
\[
\sum_{k=1}^\infty \frac{(-1)^{k-1}}{k^3 {2k \choose k} {3k \choose k}}
\left( - 2 \frac{56k^2-32k+5}{(2k-1)^2} H_{2k-1} +{\frac {280\,{k}^{
			3}-168\,{k}^{2}+33\,k-2}{2 k \left( 2\,k-1 \right) ^{2}}} \right).
\]
Notice that
\[
{\frac {280\,{k}^{3}-168\,{k}^{2}+33\,k-2}{2 k \left( 2\,k-1
		\right) ^{2}}}= \frac{5}{12} \left( \frac {6(56\,{k}^{2}-32\,k+5)}{ \left( 2\,k-1
		\right) ^{2}} \right) - \frac{1}{k}.
\]
The value for the left hand side of \eqref{eq-18} is
\[
\sum_{k=0}^\infty \left(
-4\,{\frac {H_k }{ \left( 1+k \right) ^{3}}}+
2\,{\frac {5\,k+1}{ \left( 1+k \right) ^{4}}} \right),
\]
which can be evaluated  by {\tt MultipleZetaValues}. Its value is  $-\pi^4/10 + 10 \zeta(3)$. Recall the identity \eqref{Ex21}, and we derive \eqref{eq-Chu21}.

We also  utilize Equation \eqref{eq-18} and consider the Taylor expansion at $(a,b,c,d,e)=(2,1,1,1,1)$. This time, we consider the linear combination
\[
21[a]+19[b]+\frac{115}{3}[d]+\frac{80}{3}[e].
\]
The coefficient of the right hand side of \eqref{eq-18} simplifies to the left hand of \eqref{eq-1.12}. While the coefficient of the left hand side of \eqref{eq-18} becomes
\[
\sum_{k=0}^\infty  \left(  \frac{168 H_k}{(k+1)^3}+ \frac{21}{(k+1)^4} \right),
\]
which equals $7 \pi^4/10$ by {\tt MultipleZetaValues}. This completes the proof of \eqref{eq-1.12}. \qed

\noindent {\it Proof of Theorem~\ref{Chu24}.} Once again, we utilize the identity \eqref{eq-18}.
To prove the first identity, we investigate the Taylor expansion at $(a,b,c,d,e)=(2,1,1,1,3/2)$ and the linear combination  $[a]+[b]+[d]$.
The value for the right hand side of \eqref{eq-18} is
\[
\sum_{k=1}^\infty \frac{(-1)^k}{k^2 {2k \choose k} {3k \choose k}}
\left( 8 \frac{7k-2}{2k-1} (2 H_{2k-1} - H_{k-1}) - \frac{1}{2} \frac{336k^3-160k^2+16k-1}{k (2k-1)^2} \right)	.
\]
Notice that
\begin{multline*}
- \frac{1}{2} \frac{336k^3-160k^2+16k-1}{k (2k-1)^2} \\
= \frac{1}{4} \cdot \frac{-6(56\,{k}^{2}-32\,k+5)}{k \left( 2\,k-1 \right) ^{2}}  + 2 \cdot \frac{-6(7\,k-2)}{2\,k-1} + \frac{8}{k}
\end{multline*}
In addition to \eqref{Ex21}, we also have \eqref{Ex24}.
While the value for the left hand side of \eqref{eq-18} is
\[
\sum_{k=0}^\infty  \left( {4\frac {H_k }{ \left( 1+k \right) ^{2}}}-{4\frac {H_{2k} }{ \left( 1+k \right) ^{2}}}
+{\frac {12\,{k}^{2}+12\,k+1}{ \left( 1+2\,k \right)  \left( 1+k
		\right) ^{3}}} \right).
\]
By {\tt MultipleZetaValues}, this sum reduces to $\pi^2-6 \zeta(3)^2$. Therefore, \eqref{eq-Chu24a} follows immediately.

For the second identity, we choose the linear combination
\[
[e^2] - 2 [de] - 8 [d^2] + 6 [cd].
\]
The value for the right hand side of \eqref{eq-18} is
\[
\sum_{k=1}^\infty \frac{(-1)^k}{k^2 {2k \choose k} {3k \choose k}} \left(
-8\,{\frac { \left( 7\,k-2 \right)H_{k-1}^{(2)} }{2\,k-1}}-{\frac {336\,{k}^{4}-320\,{k}^{3}+56\,{k}^{2}+19\,k-6}{{k}^{2}
		\left( 2\,k-1 \right) ^{2}}}.
\right)
\]
Notice that
\begin{multline*}
-{\frac {336\,{k}^{4}-320\,{k}^{3}+56\,{k}^{2}+19\,k-6}{{k}^{2}
		\left( 2\,k-1 \right) ^{2}}} =
 \frac{56\,{k}^{2}-32\,k+5}{k \left( 2\,k-1 \right) ^{2}} -  \frac {24(7\,k-2)}{2\,k-1} + \frac{6}{k^2}.
\end{multline*}
While the value for the left hand side of \eqref{eq-18} is
\begin{multline*}
\sum_{k=0}^\infty
\left( 36\,{\frac {H_k^{2}}{
		\left( 1+k \right) ^{2}}}-24\,{\frac { \left( (4\, {k}^{2}+6\,k+2) H_{2k} -2\,{k}^{2}-3\,k \right) H_k }{ \left( 1+2\,k \right)  \left( 1+k \right) ^{3}}}+16\,{
	\frac { H_{2k} ^{2}}{ \left(
		1+k \right) ^{2}}} \right. \\
\left.	-16\,{\frac { \left( 3+2\,k \right) kH_{2k} }{ \left( 1+2\,k \right)  \left( 1+k \right) ^{
			3}}}+4\,{\frac {k \left( 12\,{k}^{3}+36\,{k}^{2}+27\,k+4 \right) }{
		\left( 1+2\,k \right) ^{2} \left( 1+k \right) ^{4}}}\right).
\end{multline*}
By {\tt MultipleZetaValues}, this sum reduces to $2 \pi^2 - \frac{\pi^4}{90} - 4 \zeta(3)$. Therefore, \eqref{eq-Chu24b} follows immediately. \qed

\noindent
{\it Proof of Theorem~\ref{th-1.9}. } We continue to use the identity \eqref{eq-18}.
Let us consider the Taylor expansion at $(a,b,c,d,e)=(1,\frac{1}{2},\frac{1}{2},\frac{1}{2},\frac{1}{2})$ and the linear combination
\[
20[a]+12[b]+28[e].
\]
The coefficient of the right hand side of \eqref{eq-18} is
\begin{align*}
& \sum_{k=1}^\infty
\frac{(-16)^k \big( (20k^2-8k+1) H_{4k-1} - \frac{200 k^3-208 k^2+55 k-5}{2k-1} \big)}{(2k-1)^2 k^3 {2k \choose k} {4k \choose 2k}} \\
= &
\sum_{k=1}^\infty
\frac{(-16)^k \big( (20k^2-8k+1) H_{4k-1} + \frac{k(8k-5)}{2k-1} \big)}{(2k-1)^2 k^3 {2k \choose k} {4k \choose 2k}}
\\
& \qquad \qquad - 5
\sum_{k=1}^\infty
\frac{(-16)^k (20k^2-8k+1) }{(2k-1)^2 k^3 {2k \choose k} {4k \choose 2k}}.
\end{align*}
Recall the known identity \eqref{Ex1}.
While the coefficient of the left hand of \eqref{eq-18} simplifies to
\[
20 \sum_{k=0}^\infty \frac{8k+1}{(2k+1)^4},
\]
which equals $-5\pi^4/8 +70 \zeta(3)$ by {\tt MultipleZetaValues}. This completes the proof of \eqref{eq-th8-a}. 

Taking the linear combination
\[
8[a]+8[b] + 8[e],
\]
we derive \eqref{eq-th8-b} in a similar way. It is easy to see that \eqref{eq-th8-c} can be obtained by the linear combination of \eqref{eq-th8-a} and \eqref{eq-th8-b}.

For \eqref{eq-th8-d}, we consider the linear combination
\[
-16[cde]+48[de^2]-48[e^3].
\]
\qed

\noindent {\it Proof of Theorem~\ref{th-22}.} Let us recall the transformation \cite[Theorem 24]{ChZh}
\begin{align}
	& \sum_{k=0}^\infty (a+2k) \frac{(b)_k (c)_k (d)_k (e)_k}{(1+a-b)_k (1+a-c)_k (1+a-d)_k (1+a-e)_k} \nonumber \\
	= &\sum_{k=0}^{\infty} \mu_k(a ; b, c, d, e) \frac{(-1)^k}{(1+a-b)_{3 k}}  \frac{(c)_k (1+a-b-d)_k (1+a-b-e)_k}{(d+e-a)_k (1+a-d)_k (1+a-e)_k} \nonumber \\
	& \times \frac{(d)_{2k} (e)_{2k} (1+a-b-c)_{2k}}{(1+2a-b-c-d-e)_k (1+a-c)_{2k}},\label{eq-24}
\end{align}
where
{\small
	\begin{multline*}
	 \mu_k(a ; b, c, d, e)= \frac{(1 + 2a - b - c - d +3k)(a - e + k)}{1+2a - b - c - d - e + k}  \\
		 + \frac{(e +2k)(1+ a - b - c +2k)(1+ a - b - d + k)}{(1 + a - b +3k)(1+2a - b - c - d - e + k) }  \\
		 + \frac{(c+k)(d+2k)(e+2k)(1+a-b-c+2k)(1+a-b-d+k)(1+a-b-e+k)}{(d+e-a+k)(1+a-c+2k)(1+2a-b-c-d-e+k)(1+a-b+3k)_2} .
	\end{multline*}
}

Firstly, we consider the Taylor expansion at $(a,b,c,d,e)=(3/2,1,-\infty,1,1)$ and the linear combination
\[
128[e]-128[b]-128[a].
\]
The coefficient of the left hand of \eqref{eq-24} is
\begin{multline*}
\sum_{k=1}^\infty \frac{256^k \big( (22k-1) H_{2k-1} -3\,{\frac {88\,{k}^{2}-78\,k+7}{2\,k-1}}  \big)}{(2k-1) k^2 {3k \choose k}{6k \choose 3k}}
\\
= \sum_{k=1}^\infty \frac{256^k \big( (22k-1) H_{2k-1} + \frac{90k-15}{2k-1} \big)}{(2k-1) k^2 {3k \choose k}{6k \choose 3k}} 
-6  \sum_{k=1}^\infty \frac{256^k (22k-1)}{(2k-1) k^2 {3k \choose k}{6k \choose 3k}}.
\end{multline*}

Denote the series on the left hand side of \eqref{eq-24} by $\Omega(a;b,c,d,e)$. We have the transformation formula \cite[Theorem 5]{ChZh}:
\begin{multline}\label{eq-tr}
\Omega(a;b,c,d,e) = \\
\frac{a-c}{1+2a-b-c-d-e} \Omega(1+2a-b-d-e;1+a-b-e, 1+a-b-d, 1+a-d-e, c).
\end{multline}
The coefficient of the right hand of \eqref{eq-tr} is
\[
-256\,{\frac { \left( 6\,k+1 \right)  \left( -1 \right) ^{k}}{ \left( 
		2\,k+1 \right) ^{3}}},
\]
which equals $-768G  +16 \pi^3$ by {\tt MultipleZetaValues}. This completes the proof of \eqref{eqt-1.9a}. 

In a similar way, we can prove \eqref{eqt-1.9b} by the combination
\[
256[e]-256[b]-384[a],
\]
and prove \eqref{eqt-1.9d} by the combination
\[
-256[e^2]+128 [de]. \tag*{\qed}
\]

\noindent {\it Proof of Theorem~\ref{th-1.10}.}
We use \eqref{eq-24} and fix the Taylor expansion at $(a,b,c,d,e)=(\frac{3}{2}, 1, \frac{1}{2}, 1, 1)$.
The first combination we consider is
\[
-128([a]+[b]+[c]+2[e]).
\]
The coefficient of the right hand side of \eqref{eq-24} is the same as the left hand side of \eqref{eq-4.33}. By aid of the transform \eqref{eq-tr} once again, we find that the coefficient of  the left hand side of \eqref{eq-24}
is
\[
-768\sum_{k=0}^\infty \frac{1}{(2k+1)^4}.
\]
By {\tt MultipleZetaValues}, it equals $-8\pi^4$. This completes the proof of \eqref{eq-4.33}.

The second combination we consider is
\[
-128(4[a]+5[b]+4[c]+7[e]).
\]
By a similar argument as before, we find that the series on the left hand side of \eqref{eq-4.34} is equal to
\[
-3072\sum_{k=0}^\infty \frac{1}{(2k+1)^4},
\]
which completes the proof of \eqref{eq-4.34}.

The third combination we consider is
\[
128([a^2]+[ab]+[ac]+[ae]+[b^2]+[bc]+[be]+[c^2]+[ce]+[e^2]).
\]
By a similar argument as before, we find that the series on the left hand side of \eqref{eq-4.35} is equal to
\[
1024\sum_{k=0}^\infty \frac{1}{(2k+1)^5},
\]
which equals $992\zeta(5)$, completing the proof of \eqref{eq-4.35}.

Finally, we consider the combination
\begin{multline*}
-128 \big( [a^3]+[a^2b]+[a^2c]+[a^2e]+[ab^2]+[abc]+[abe]+[ac^2]+[ace]+[ade]-[ae^2]+[b^3]\\
+[b^2c]+[b^2e]+[bc^2]+[bce]+[bde]-[be^2]+[c^3]+[c^2e]+[cde]-[ce^2]+[de^2]-2[e^3] \big).
\end{multline*}
By a similar argument as before, we find that the series on the left hand side of \eqref{eq-4.36} is equal to
\[
-512 \sum_{k=0}^\infty \frac{8H^{(3)}_{2k}-H^{(3)}_{k}}{(2k+1)^3} - 2048 \sum_{k=0}^\infty \frac{1}{(2k+1)^6},
\]
which equals $-1568 \zeta(3)^2$, completing the proof of \eqref{eq-4.36}. \qed

\noindent {\it Proof of Theorem~\ref{th-1.11}.}
We utilize the transformation \eqref{eq-18} and consider the Taylor expansion at $(a,b,c,d,e)=(\frac{3}{2}, 1, \frac{1}{2}, \frac{1}{2}, 1)$.

The first combination is
\[
8[a]+8[b]+6[d]+11[e].
\]
We find that the series on the left hand side of \eqref{eq-1.11a} equals
\[
\sum_{k=0}^\infty \left(
\frac{9(4k+3)H_{2k}}{(2k+1)^2(k+1)^2}
+   \frac{16k^3+24k^2+21k+8}{(2k+1)^3(k+1)^3} \right).
\]
By {\tt MultipleZetaValue}, it equals $11\zeta(3)$, completing the proof of \eqref{eq-1.11a}.

The second combination is
\[
36[a]+40[b]+28[d]+46[e].
\]
We find that the series on the left hand side of \eqref{eq-1.11b} equals
\[
\sum_{k=0}^\infty \left( 
 \frac{42(4k+3)H_{2k}}{(2k+1)^2(k+1)^2}
+  \frac{2(32k^3+48k^2+45k+18)}{(2k+1)^3(k+1)^3} \right).
\]
By {\tt MultipleZetaValue}, it equals $50\zeta(3)$, completing the proof of \eqref{eq-1.11b}.
\qed

\noindent{\it Proof of Theorem~\ref{th-1.12}.}
We utilize the transformation \eqref{eq-18} and consider the Taylor expansion at $(a,b,c,d,e)=(1, 1/2, 1, 1/2, 1/2)$. 

The first combination is 
\[
-[c]+[d] -[e].
\]
The coefficient of the right hand side of \eqref{eq-18} is exactly the same as the left hand side of \eqref{eq-1.12a}. While coefficient of the left hand sider of \eqref{eq-18} is
\[
-2 \sum_{k=0}^\infty \frac{H_k}{(2k+1)^2},
\]
which equals $\frac{\pi^2}{2} \log 2- \frac{7}{2} \zeta(3)$ by {\tt MultipleZetaValue}. This proves \eqref{eq-1.12a}.

The second combination is 
\[
-2[b]+2[c]
\]
The coefficient of the right hand side of \eqref{eq-18} is
\begin{multline*}
	\sum_{k=1}^\infty \frac{(-16)^{k} \big( P(k) (2H_{6k-1}-H_{3k-1}-3H_{k-1}) -224k^3+168k^2-24k-5/2 \big)}{k^2(2k-1)(4k-1)(4k-3) \bi{3k}k\bi{6k}{3k}} \\
	= \sum_{k=1}^\infty \frac{(-16)^{k} \big( P(k) (2H_{6k-1}-H_{3k-1}-3H_{k-1}) -64k^2+46k -17/2 \big)}{k^2(2k-1)(4k-1)(4k-3) \bi{3k}k\bi{6k}{3k}} \\
	-2 \sum_{k=1}^\infty \frac{(-16)^{k} P(k) }{k^2(2k-1)(4k-1)(4k-3) \bi{3k}k\bi{6k}{3k}}	
\end{multline*}
While coefficient of the left hand sider of \eqref{eq-18} is
\[
\sum_{k=0}^\infty \left( \frac{8(H_k-H_{2k})}{(1+2k)^2} + \frac{8k}{(1+2k)^3} \right),
\]
which equals $\pi^2/2 - 2 \pi^2 \log(2) + 7 \zeta(3)$ by {\tt MultipleZetaValue}. Noting the formula \eqref{Ex25}, we derive \eqref{eq-1.12b}. \qed

\noindent {\it Proof of Theorem~\ref{th-1.13}.} 
We utilize \eqref{eq-14} and consider the Taylor expansion at $(a,b,c,d,e)=(3/2, 1, 1, 1, 1/2)$.

The first combination is
\[
[a]+[c]+2[d]+[e].
\]
The coefficient of the right hand side of \eqref{eq-14} is 
\[
\sum_{k=1}^\infty (a_{2k-1}+a_{2k}) = \sum_{k=1}^\infty a_k,
\]
where $a_k$ is the summand of the left hand side of \eqref{eq-1.13a}. By the transformation \eqref{eq-tr}, the coefficient of the left hand side of \eqref{eq-14} is
\[
\sum_{k=0}^\infty \frac{6}{(2k+1)^4} = \frac{\pi^4}{16},
\]
completing the proof of \eqref{eq-1.13a}.

The second combination is
\[
[a^2] +[ac] + [ad]+ [ae]+[c^2]+[cd]+[ce]+[d^2]+[de]+[e^2].
\]
The coefficient of the right hand side of \eqref{eq-14} is 
\[
\sum_{k=1}^\infty (a_{2k-1}+a_{2k}) = \sum_{k=1}^\infty a_k,
\]
where $a_k$ is the summand of the left hand side of \eqref{eq-1.13b}. By the transformation \eqref{eq-tr}, the coefficient of the left hand side of \eqref{eq-14} is
\[
\sum_{k=0}^\infty \frac{8}{(2k+1)^5} = \frac{31}{4} \zeta(5),
\]
completing the proof of \eqref{eq-1.13b}.

The third combination is
\begin{multline*}
[a^3]+[a^2c]+[a^2d]+[a^2e]+[abd]+[ac^2]+[acd]+[ace]-[ad^2]+[ade]+[ae^2]+[bcd]\\
+[bd^2]+[bde]+[c^3]+[c^2d]+[c^2e]-[cd^2]+[cde]+[ce^2]-2[d^3]-[d^2e]+[de^2]+[e^3].
\end{multline*}
The coefficient of the right hand side of \eqref{eq-14} is 
\[
\sum_{k=1}^\infty (a_{2k-1}+a_{2k}) = \sum_{k=1}^\infty a_k,
\]
where $a_k$ is the summand of the left hand side of \eqref{eq-1.13c}. By the transformation \eqref{eq-tr}, the coefficient of the left hand side of \eqref{eq-14} is
\[
\sum_{k=0}^\infty \left( \frac{4(8H^{(3)}_{2k}-H^{(3)}_k)}{(2k+1)^3} + \frac{16}{(2k+1)^6} \right)= \frac{49}{4} \zeta(3)^2,
\]
completing the proof of \eqref{eq-1.13c}.

\noindent {\it Proof of Theorem~\ref{th-1.14}.} 
We utilize \eqref{eq-24} and consider the Taylor expansion at $(a,b,c,d,e)=(1,1/2,-\infty,1/2,1)$.

The first combination is
\[
-8([a]+[b]).
\]
The coefficient of the right hand side of \eqref{eq-24} is 
\begin{align*}
& \sum_{k=1}^\infty \frac{16^k {4k \choose 2k} \big( Q(k) H_{k-1} - (88 k^3 -136 k^2+68 k-9)/(2k-1) \big) }{k (4k-1)(4k-3){3k \choose k} {6k \choose 3k}} \\
= & \sum_{k=1}^\infty \frac{16^k {4k \choose 2k} \big( Q(k) H_{k-1} + (6k-1)(4k-3)/(2k-1) \big) }{k (4k-1)(4k-3){3k \choose k} {6k \choose 3k}} \\
& - 2  \sum_{k=1}^\infty \frac{16^k {4k \choose 2k} Q(k) }{k (4k-1)(4k-3){3k \choose k} {6k \choose 3k}} \\
= &  \sum_{k=1}^\infty \frac{16^k {4k \choose 2k} \big( Q(k) H_{k-1} + (6k-1)(4k-3)/(2k-1) \big) }{k (4k-1)(4k-3){3k \choose k} {6k \choose 3k}} - 4 \pi.
\end{align*}
The coefficient of the left hand side of \eqref{eq-24} is
\[
\sum_{k=0}^\infty \left( \frac{8 (-1)^k H_k}{2k+1} - \frac{8 (-1)^k (4 k+1)}{(2k+1)^2} \right)  = 16 G - 4 \pi - 4 \pi \log 2,
\]
completing the proof of \eqref{eq-1.14a}.

The second combination is
\[
-8 ( [a]+[b]+[d]/3).
\]
The coefficient of the right hand side of \eqref{eq-24} is 
\begin{align*}
	& \sum_{k=1}^\infty \frac{16^k {4k \choose 2k} \big( Q(k) H_{2k-1} - (264 k^3-342 k^2+157 k-21)/(6k-3) \big) }{k (4k-1)(4k-3){3k \choose k} {6k \choose 3k}} \\
	= & \sum_{k=1}^\infty \frac{16^k {4k \choose 2k} \big( Q(k) H_{2k-1} + (6k-1)(k-3)/(6k-3) \big) }{k (4k-1)(4k-3){3k \choose k} {6k \choose 3k}} \\
	& - 2  \sum_{k=1}^\infty \frac{16^k {4k \choose 2k} Q(k) }{k (4k-1)(4k-3){3k \choose k} {6k \choose 3k}} \\
	= &  \sum_{k=1}^\infty \frac{16^k {4k \choose 2k} \big( Q(k) H_{2k-1} + (6k-1)(4k-3)/(2k-1) \big) }{k (4k-1)(4k-3){3k \choose k} {6k \choose 3k}} - 4 \pi.
\end{align*}
The coefficient of the left hand side of \eqref{eq-24} is
\[
\sum_{k=0}^\infty \left( \frac{8 (-1)^k H_k}{3(2 k+1)}  - \frac{8 (-1)^k (4 k+1)}{(2k+1)^2} \right) = \frac{32}{4} G - 4 \pi - \frac{4}{3} \pi \log 2.
\]
completing the proof of \eqref{eq-1.14b}.

The third combination is
\[
-8 (2[a]+[b]+2[d]/3).
\]
The coefficient of the right hand side of \eqref{eq-24} is
{\small 
\begin{align*}
	& \sum_{k=1}^\infty \frac{16^k {4k \choose 2k} \big( Q(k) (2H_{6k-1}-H_{3k-1}) - 2(396 k^3-540 k^2+247 k-33)/(6k-3) \big) }{k (4k-1)(4k-3){3k \choose k} {6k \choose 3k}} \\
	= & \sum_{k=1}^\infty \frac{16^k {4k \choose 2k} \big( Q(k) (2H_{6k-1}-H_{3k-1}) + 4(18 k^2- 20 k+3)/(6k-3) \big) }{k (4k-1)(4k-3){3k \choose k} {6k \choose 3k}} \\
	& - 6  \sum_{k=1}^\infty \frac{16^k {4k \choose 2k} Q(k) }{k (4k-1)(4k-3){3k \choose k} {6k \choose 3k}} \\
	= &  \sum_{k=1}^\infty \frac{16^k {4k \choose 2k} \big( Q(k) (2H_{6k-1}-H_{3k-1}) + (6k-1)(4k-3)/(2k-1) \big) }{k (4k-1)(4k-3){3k \choose k} {6k \choose 3k}} - 12 \pi.
\end{align*}}\\
The coefficient of the left hand side of \eqref{eq-24} is
\[
\sum_{k=0}^\infty \left( \frac{(-1)^k (32 H_{2k} - 32 H_k/3)}{2k+1} - \frac{16 (-1)^k (1+6 k)}{(2 k+1)^2} \right) = \frac{64}{4} G - 12 \pi + \frac{4}{3} \pi \log 2.
\]
completing the proof of \eqref{eq-1.14c}. \qed

\noindent {\it Proof of Theorem~\ref{th-1.15}.} 
We utilize \eqref{eq-14} and consider the Taylor expansion at $(a,b,c,d,e)=(3/2,1/2,1,1/2,1)$.
Since  the proof is similar to that of Theorem~\ref{th-1.14}, we only list the combinations.

For \eqref{eq-1.15a}, the combination is
\[
4[a]+4[d]+10[e].
\]
For \eqref{eq-1.15b}, the combination is
\[
3[a] +2 [d] + 7[e].
\]
For \eqref{eq-1.15c}, the combination is
\[
2[bd]-\frac{1}{2}[ce]-4[d^2]+[e^2]
\]
For \eqref{eq-1.15d}, the combination is
\[
2[a^2]+2[ad]+4[ae]+4[bd]+[ce]-6[d^2]+4[de]+6[e^2].
\tag*{\qed}
\]

\noindent {\it Proof of Theorem~\ref{th-16}.} 
We utilize Theorem~21 of \cite{ChZh} and consider the Taylor expansion at $(a,b,c,d,e)=(1,1/2,1/2,1/2,1)$.
Since  the proof is similar to that of Theorem~\ref{th-1.14}, we only list the combinations.

For \eqref{eq-16a}, the combination is
\[
192([c]-2[e]).
\]
For \eqref{eq-16b}, the combination is
\[
-384(5[a]+3[c]+4[e]).
\tag*{\qed}
\]

\noindent {\it Proof of Theorem~\ref{th-17}.} 
We utilize Theorem~27 of \cite{ChZh} and consider the Taylor expansion at $(a,b,c,d,e)=(2, 1, 1, 1, 3/2)$. 
The proof is similar to that of Theorem~\ref{th-1.14} except that there is one extra term consisting of the product of gamma functions on the left hand side of Theorem~27. Once again, we only list the combinations.

For \eqref{eq-17a}, the combination is
\[
15[d]+3[a].
\]
For \eqref{eq-17b}, the combination is
\[
-291[d^2]+69[cd]-51[ad]-17[a^2].
\]
For \eqref{eq-17c}, the combination is
\[
9[d^3]-9[cd^2]+3[bcd].
\tag*{\qed}
\]

\noindent {\it Proof of Theorem~\ref{th-18}.} 
We utilize Theorem~9 of \cite{ChZh} and consider the Taylor expansion at $(a,b,c,d,e)=(2, 2, 1/2, 1/2, 1/2)$. 
The proof is similar to that of Theorem~\ref{th-1.14}, except that additional effort is required for the evaluation.

The combination is given by
\[
-2[e]- \frac{10}{9} [b]- \frac{26}{27}[a].
\]
The right hand side is exactly \eqref{eqt-18}, while the left hand side equals
\begin{multline}
\label{eq18r}
\sum_{k=0}^\infty 
\Big( -8\,{\frac { \left( k+1 \right) ^{2}H_k }{
		\left( 2\,k+1 \right) ^{3} \left( 2\,k+3 \right) ^{3}}}-120\,{\frac {
		\left( k+1 \right) ^{2}H_{2k}}{ \left( 2
		\,k+1 \right) ^{3} \left( 2\,k+3 \right) ^{3}}} \\
	-2\,{\frac { \left( k+1
		\right)  \left( 136\,{k}^{3}+388\,{k}^{2}+334\,k+39 \right) }{
		\left( 2\,k+3 \right) ^{4} \left( 2\,k+1 \right) ^{4}}}
\Big).
\end{multline}
The package  {\tt MultipleZetaValue} cannot treat the above sum. We use the following identity
\[
\sum_{k=0}^\infty (a_{k+1}-a_k) H_k = - \sum_{k=1}^\infty \frac{a_k}{k}.
\] 
For the first summand, we have
\[
-8\,{\frac { \left( k+1 \right) ^{2}}{ \left( 2\,k+1 \right) ^{3}
		\left( 2\,k+3 \right) ^{3}}}
= \Delta \left(-\,{\frac {k \left( 2\,k+3 \right) }{ 4\left( 2\,k+1 \right) ^{3}}} \right) - \frac{1}{4 \left( 2\,k+1 \right) ^{2}},
\]
so that
\[
\sum_{k=0}^\infty -8\,{\frac {  \left( k+1 \right) ^{2} H_k}{ \left( 2\,k+1 \right) ^{3}
		\left( 2\,k+3 \right) ^{3}}}
= \sum_{k=1}^\infty \frac  {\left( 2\,k+3 \right) }{ 4\left( 2\,k+1 \right) ^{3}} 
- \sum_{k=0}^\infty  \frac{H_k}{4 \left( 2\,k+1 \right) ^{2}}.
\]
Both of the two sums on the right hand side can be evaluated by {\tt Mathematics}, which leads to 
\[
\frac{1}{32}  (-24 + \pi^2 (1 +  \log 4)).
\]
Similarly,
\begin{align*}
& \sum_{k=0}^\infty \left( -120\,{\frac { \left( k+1 \right) ^{2} H_{2k}}{ \left( 2\,k+1 \right) ^{3}
		\left( 2\,k+3 \right) ^{3}}} \right) \\
& = \sum_{k=1}^\infty \frac{15}{4} \frac{k (2k+3)}{(2k+1)^3} \left(\frac{1}{2k}+\frac{1}{2k-1} \right)
- \sum_{k=0}^\infty \frac{15 H_{2k}}{4 (2 k+1)^2} 
\\
&= \frac{15}{64} ( \pi^2 -28+ 14\zeta(3)),
\end{align*}
and
\begin{align*}
& -2\sum_{k=0}^\infty
{\frac { \left( k+1 \right)  \left( 136\,{k}^{3}+388\,{k}^{2}+334
		\,k+39 \right) }{ \left( 2\,k+3 \right) ^{4} \left( 2\,k+1 \right) ^{4
}}} 
 = \frac{1}{64} (420 - 17 \pi^2 - 266 \zeta(3)).
\end{align*}
Summing up these three values, we arrive at \eqref{eqt-18}. \qed

\noindent {\it Proof of Theorem~\ref{th-19}.} 
We utilize \eqref{eq-18} and consider the Taylor expansion at $(a,b,c,d,e)=(1, 1/2, 1/2, 1/2, -\infty)$. 
The proof is similar to that of Theorem~\ref{th-1.14}. We only give the combination:
\[
-64[c]-32[b]-64[a]. \tag*{\qed}
\]

\noindent {\it Proof of Theorem~\ref{th-20}.} 
We utilize \eqref{eq-14} and consider the Taylor expansion at $(a,b,c,d,e)=(1, 1/2, 1/2, 1/2, -\infty)$. 
The proof is similar to that of Theorem~\ref{th-1.14}. We only give the combinations:
\[
64[d]+32[c]+64[a], \quad 128[d]-32[c]+64[a],
\quad -128[d^2]+64[bd]. \tag*{\qed}
\]

\noindent {\it Proof of Theorem~\ref{th-21}.} 
We utilize \cite[Theorem~31]{ChZh}  and consider the Taylor expansion at $(a,b,c,d,e)=(2,1,1,1,3/2)$. 
The proof is similar to that of Theorem~\ref{th-1.14}. We only give the combinations:
\[
\frac{1}{2}  ( [a]+[d]+[e] ), \quad -\frac{1}{4} ([a]+5[d]),
\]
\[
3[d^2]+3[cd]+3[ad]+[a^2], \quad 25[d^2]+[cd]+9[ad]+3[a^2]. \tag*{\qed}
\]

\noindent {\it Proof of Theorem~\ref{th-22n}.} 
We utilize \eqref{eq-14} and consider the Taylor expansion at $(a,b,c,d,e)=(1,1/2,1/2,1/2,1/2)$. 
The proof is similar to that of Theorem~\ref{th-1.14}. We only give the combinations:
\[
-64  ( 2[d]-[a] ), \quad -64 ([a]+[d]+[e]),
\]
\[
32([ce]-2[e^2]), \quad 16(-[bd]+[ce]+2[d^2]-2[e^2]). \tag*{\qed}
\]

\noindent {\it Proof of Theorem~\ref{th-23}.} 
We consider the Taylor expansion at $(a,b,c,d,e)=(2,1,1,1,1)$ and the combinations:
\[
2[a]-2[b]+10[e], \quad 2[a]+8[e], \quad 1514[e^2]+134 [de] + 594 [ae] +198 [a^2].
\]

\end{document}